\documentclass[a4paper,12pt]{amsart}
\usepackage{amssymb}
\usepackage{amsmath}
\usepackage{color}

\makeatletter
\@namedef{subjclassname@2020}{%
  \textup{2020} Mathematics Subject Classification}
\makeatother

\usepackage{mathtools}
\usepackage{stmaryrd}

\newcommand{\ignore}[1]{}

\newcommand{\Z}{\mathbb{Z}}
\newcommand{\Q}{\mathbb{Q}}

\newcommand{\F}{\mathbb{F}}

\newcommand{\ds}{\displaystyle}

\DeclareMathOperator*\Res{Res}

\def\ds{\displaystyle}
\def\Li{{\rm Li}}
\def\LL{\pounds}

\newcommand{\leg}[2]{\left(\frac{#1}{#2}\right)}

\newtheorem{dummy}{Dummy}

\newtheorem{lemma}[dummy]{Lemma}
\newtheorem{theorem}[dummy]{Theorem}
\newtheorem*{theorem*}{Theorem}

\newtheorem{cor}[dummy]{Corollary}

\theoremstyle{definition}

\theoremstyle{remark}
\newtheorem{rem}[dummy]{Remark}

\newtheorem*{rem*}{Remark}

\begin{document}
\bibliographystyle{amsalpha}

\author{S.~Mattarei}
\email{sandro.mattarei@unimib.it}
\urladdr{https://www.unimib.it/sandro-mattarei}
\address{Dipartimento di Matematica e Applicazioni\\
  Universit\`a di Milano--Bicocca\\
  via Roberto Cozzi\\
  20125 Milano\
  Italy}

  \author{R.~Tauraso}
\email{tauraso@mat.uniroma2.it}
\urladdr{https://www.mat.uniroma2.it/\~{ }tauraso/}
\address{Dipartimento di Matematica\\
  Universit\`a di Roma ``Tor Vergata'' \\
  via della Ricerca Scientifica\\
  00133 Roma\\
  Italy}

\title{Congruences for sums involving $\binom{rk}{k}$}

\begin{abstract}
We primarily investigate congruences modulo $p$ for finite sums of the form
$\sum_k\binom{rk}{k}x^k/k$
over the ranges $0<k<p$ and $0<k<p/r$,
where $p$ is a prime larger than the positive integer $r$.
Here $x$ is an indeterminate, thus allowing specialization to numerical congruences where $x$ takes certain algebraic numbers as values.
We employ two different approaches that have complementary strengths.
In particular, we obtain congruences modulo $p^2$ for the sum
$\sum_{0<k<p}\binom{rk}{k}x^k$,
expressed in terms of finite polylogarithms of certain quantities related to $x$.

\end{abstract}


\maketitle

\thispagestyle{empty}

\section{Introduction}

This may be viewed as the third of a series of papers in which the authors develop polynomial congruences,
modulo a prime or a power of a prime, where the polynomial coefficients are closely related to binomial coefficients.
As a simple, paradigmatic instance, one may take the congruence
\begin{equation}\label{eq:central_pol}
\sum_{0\le k<q/2}\binom{2k}{k}x^k\equiv(1-4x)^{(q-1)/2}\pmod{p},
\end{equation}
which is valid for any power $q$ of an odd prime $p$.
Although Equation~\eqref{eq:central_pol} admits various easy direct proofs, it seems natural to deduce it from the analogous
power-series identity
%\begin{equation}\label{eq:central_gf}
$
\sum_{k=0}^\infty\binom{2k}{k}x^k
=
(1-4x)^{-1/2}.
% 1/\sqrt{1-4x}.
$
A procedure for this deduction,
which may be called {\em truncation and reduction modulo $p$} of a power-series identity,
was extensively investigated in~\cite{MatTau:truncation} for variations of binomial coefficients of the form $\binom{2k}{k}$.

It is worth noting that quite a number of papers
were devoted to prove congruences closely related to Equation~\eqref{eq:central_pol}, sometimes with considerable effort,
that were numerical rather than polynomial, meaning that they were stated for particular rational values of $x$.
A number of those earlier results by various authors can be found cited in the survey~\cite{Sun:sums_higher_Catalan}.
An advantage of polynomial congruences that was emphasized in~\cite{MatTau:truncation}
is that they can be differentiated, and sometimes integrated explicitly, two processes that are unavailable with numerical congruences.
Subsequent substitution of rational,
or more generally algebraic numbers, for the indeterminate $x$
yields numerical congruences such as
$
\sum_{0\le k<p}\binom{2k}{k}k^{-3}\equiv 2B_{p-3}/3\pmod{p},
$
where $p>3$ is a prime and $B_{p-3}$ is a Bernoulli number
(see~~\cite[Eq.~(42)]{MatTau:truncation}),
which would be quite challenging to prove directly.

The investigation in~\cite{MatTau:truncation} was extended in~\cite{MatTau:truncation-cubic}
to series involving binomial coefficients of the form $\binom{3k}{k}$.
Thus, the authors produced polynomial congruences, modulo a prime $p$,
for certain partial sums of the generating series of $\binom{3k}{k}$ and, more generally, of
the power series $\sum_{k=0}^{\infty}\binom{3k+e}{k}x^k$, where $e\ge 0$ is an integer.
For example, as a special case of~\cite[Corollary~3.4]{MatTau:truncation-cubic}, we showed
that the congruence
\begin{equation*}
2(2+c)
(1-c)^{q-1}
\sum_{0\le k<q/3}\binom{3k}{k}x^k
\equiv
1+3(1-4c)^{(q-1)/2}
+2c^q
\pmod{p}
\end{equation*}
holds in the polynomial ring $\Z[c]$, where $x=c^2/(1-c)^3$,
and $q$ is a power of an odd prime $p$.
This should be viewed as a polynomial-congruence version of the identity
\[
\sum_{k=0}^{\infty}\binom{3k}{k}x^k
=
\frac{1-c}{2(2+c)}\left(1+3(1-4c)^{-1/2}\right)
\]
in the ring of formal power series $\Q[[x]]$.
%The method for obtaining a polynomial congruence from a power series identity drew on ideas introduced
%in~\cite{MatTau:truncation} in the simpler setting of sums that involve the
%central binomial coefficients $\binom{2k}{k}$.
In this case as well, our polynomial approach demonstrates its strength
by allowing one to recover certain numerical congruences in the literature, and to produce more, simply by specializing $x$ to certain numbers.

The above introduction of an accessory indeterminate $c$ related to $x$
resulted from
the need to be able to solve a cubic equation
satisfied by the series
$\sum_{k=0}^{\infty}\binom{3k}{k}x^k$.
More generally, for any positive integer $r$ the power series $y=y_r(x)=\sum_{k=0}^\infty\binom{rk}{k}x^k$ satisfies
\begin{equation*}\label{eq:y_equation}
(y-1)\bigl((r-1)y+1\bigr)^{r-1}-r^rxy^r=0.
\end{equation*}
As explained in~\cite[Section~2]{MatTau:truncation-cubic}, when we substitute
$x=-c^{r-1}/(c-1)^r$
the left-hand side of the resulting equation
has a linear factor
$y-(c-1)/(c-1+r)$.
Then the series $y_r(x)$ of our interest, now expressed in terms of $c$, is a root of the other factor, of degree $r-1$.
This decrease in degree is what allowed the case $r=3$ to be dealt with effectively in~\cite{MatTau:truncation-cubic}.
One of the goals of the present paper is to explore ways in which the general case of arbitrary $r$ can be treated.

Before pursuing those, in Section~\ref{sec:3ksuk_series} we discuss how to obtain congruences modulo $p$ for the polynomials
$\sum_{0<k<p}\binom{3k}{k}x^k/k$
and
$\sum_{0<k<p/3}\binom{3k}{k}x^k/k$,
which are Equations~\eqref{eq:3kksuk-pol-long} and~\eqref{eq:3kksuk-pol-short}.
One way is to apply the method of {\em truncation and reduction modulo $p$} developed in~\cite{MatTau:truncation} and~\cite{MatTau:truncation-cubic}
to a primitive of the power series $(y_r(x)-1)/x$, which we produce in Theorem~\ref{thm:rkksuk-series}.
Not unexpectedly, the resulting congruences involve instances of the {\em truncated logarithm} polynomial $\LL_1(t)=\sum_{0<k<p}t^k/k$.

However, the core of this paper focuses on two new and more powerful methods to access such and further congruences,
with the discussion of Section~\ref{sec:3ksuk_series} playing more of an inspirational role than technical service.
Key to both methods is the realization that one should make simultaneous use of {\em all} the algebraic functions $c=c_1,\ldots,c_r$ of $x$ defined by the equation
$x(c-1)^r+c^{r-1}=0$,
associated to the substitution discussed earlier.
As an example, in Theorem~\ref{thm:rkksuk} we find
\[
\sum_{0<k<p} \binom{rk}{k} \frac{x^k}{k}
\equiv
-rx^p\sum_{i=1}^{r}
\LL_1(c_i) \pmod{p}
\]
if $p>r$.
We prove further results of this type in Section~\ref{sec:pounds} based on our first new method, which we introduce in Section~\ref{sec:symmetric}.
As an application of one of them, we find that if $p>r$, and $a$ is an integer such that
the polynomial $a(c-1)^r+c^{r-1}\in\F_p[c]$ has $r$ distinct roots
in the field $\F_p$ of $p$ elements, then
$\sum_{0<k<p}\binom{rk}{k}a^k\equiv 0\pmod{p}$.

Our second new method, which we present in Section~\ref{sec:mod_p^2}, has an overlapping scope with the first.
Where both methods apply, the second produces stronger congruences, such as a congruence modulo $p^2$ for the sum
$\sum_{0<k<p} \binom{rk}{k}x^k/k$ considered above, see Theorem~\ref{thm:rkksukmod2}.
Similar congruences modulo $p^2$ for the sum $\sum_{0<k<p} \binom{rk}{k}x^k$ follow by differentiation,
thus strengthening some of the congruences obtained in~\cite{MatTau:truncation-cubic}.

To illustrate the progress enabled by our new approaches, in the final Section~\ref{sec:numerical} we present a sample of numerical congruences
obtained from our algebraic congruences upon specialization of the indeterminate $x$ to various rational numbers.

\section{The power series $\sum_{k=1}^{\infty}\binom{rk}{k}x^k/k$}
\label{sec:3ksuk_series}

The power series $y_r(x)=\sum_{k=0}^\infty\binom{rk}{k}x^k$
mentioned in the Introduction can be conveniently expressed
in terms of the more fundamental power series
\begin{equation}\label{eq:B_r_def}
\mathcal{B}_r(x)
=\sum_{k=0}^{\infty}\frac{1}{rk+1}\binom{rk+1}{k}x^k
=\sum_{k=0}^{\infty}\frac{1}{(r-1)k+1}\binom{rk}{k}x^k,
\end{equation}
see~\cite[Equation~(5.58)]{GKP},
which satisfies the simpler equation~\cite[Example~6.2.6]{Stanley:EC2}
\begin{equation}\label{eq:B_r_equation}
\mathcal{B}_r(x)=1+x\mathcal{B}_r(x)^r.
\end{equation}
Hence the power series $\mathcal{B}_r(x)$ is one solution of the equation
$z=1+xz^r$ for the unknown $z$;
in fact, that is the only solution in the ring
of formal power series $\Q[[x]]$.
In this paper we will find it useful to simultaneously  consider all the $r$ solutions of this equation,
in an appropriate setting that will be explained in Section~\ref{sec:symmetric}.

Now we have
\[
y_r(x)=\sum_{k=0}^\infty\binom{rk}{k}x^k
=
\frac{1}{1-r+r\mathcal{B}_r(x)^{-1}}
=
\frac{\mathcal{B}_r(x)}{r-(r-1)\mathcal{B}_r(x)},
\]
see~\cite[Section~2]{MatTau:truncation-cubic} for details of this derivation.
The following result shows how $\mathcal{B}_r(x)$ also plays a crucial role
in expressing a primitive of the power series $(y_r(x)-1)/x$.

\begin{theorem}\label{thm:rkksuk-series}
For any positive integer $r$, in the formal power series ring $\Q[[x]]$ we have
\[
\sum_{k=1}^{\infty} \binom{rk}{k} \frac{x^k}{k}=r\log(\mathcal{B}_r(x)).
\]
\end{theorem}

\begin{proof}
Since both power series have no constant term, it suffices to show
that the derivative of the right-hand side equals
$\sum_{k=1}^\infty\binom{rk}{k}x^{k-1}$.
Equivalently, but slightly more conveniently, we apply the differential operator
$x\cdot d/dx$ to the right-hand side.
We find
\[
x\,\frac{d}{dx}r\log(\mathcal{B}_r(x))
=
\frac{rx}{\mathcal{B}_r(x)}\cdot
\frac{d}{dx}\mathcal{B}_r(x)
=
\frac{rx\mathcal{B}_r(x)^{r-1}}{1-rx\mathcal{B}_r(x)^{r-1}}
=
-1+\frac{1}{1-rx\mathcal{B}_r(x)^{r-1}},
\]
where the derivative of $\mathcal{B}_r(x)$ was computed by implicit differentiation of
Equation~\eqref{eq:B_r_def}.
Now using Equation~\eqref{eq:B_r_def} we conclude
\[
x\,\frac{d}{dx}r\log(\mathcal{B}_r(x))
=
-1+\frac{\mathcal{B}_r(x)}{\mathcal{B}_r(x)-r(\mathcal{B}_r(x)-1)}
=\sum_{k=1}^\infty\binom{rk}{k}x^k,
\]
as desired.
\end{proof}

In the special case where $r=2$, Equation~\eqref{eq:B_r_equation}
can be explicitly solved and yields
$\mathcal{B}_2(x)=\bigl(1-\sqrt{1-4x}\bigr)/(2x)$,
the generating function of the Catalan numbers.
This special case of Theorem~\ref{thm:rkksuk-series} is well known, see~\cite[Equation~(6)]{Lehmer:central_binomial}.
To facilitate comparison with results in~\cite{MatTau:truncation} and~\cite{MatTau:truncation-cubic}
we may express this special case as an identity in the ring of formal power series $\Q[[c]]$,
where $c$ is an auxiliary indeterminate $c$ related to $x$ by $x=-c/(c-1)^2$:
\begin{equation}\label{thm:2kksuk-series}
\sum_{k=1}^{\infty} \binom{2k}{k} \frac{x^k}{k}
=2\log(1-c)=2\,\Li_1\left(\frac{c}{c-1}\right).
\end{equation}
Here $\Li_1(t)=\sum_{k=1}^{\infty}t^k/k=-\log(1-t)$.

In~\cite[Equation~(32)]{MatTau:truncation} we proved a matching congruence for a finite sum,
\begin{equation}\label{thm:2kksuk-pol}
(c-1)^p\sum_{0<k<p/2}\binom{2k}{k}\frac{x^k}{k}
\equiv
2\LL_1\left(c\right)
\pmod{p},
\end{equation}
where $p$ is an odd prime.
Here $\LL_1(t)=\sum_{k=1}^{p-1}t^k/k$
is the {\em truncated logarithm}, a special instance of a
{\em finite polylogarithm} $\LL_d(t)=\sum_{k=1}^{p-1}t^k/k^d$,
where $d$ is any integer.
Note that the summation in Equation~\eqref{thm:2kksuk-pol} may equivalently be written over the range $0<k<p$,
because $\binom{2k}{k}$ is a multiple of $p$ if $p/2<k<p$.

Equation~\eqref{thm:2kksuk-pol} takes place in the polynomial ring $\Z_p[c]$,
where $\Z_p$ is the ring of $p$-adic integers and $x=-c/(c-1)^2$.
Being a congruence between polynomials was a useful feature
in~\cite{MatTau:truncation}
because degree considerations were important in the proofs.
However, the similarity with Equation~\eqref{thm:2kksuk-series}
is greatly enhanced if we divide by $(c-1)^p$ and thus rewrite Equation~\eqref{thm:2kksuk-pol}
in the equivalent form
\begin{equation}\label{thm:2kksuk-pol*}
\sum_{0<k<p/2}\binom{2k}{k}\frac{x^k}{k}
\equiv
2\LL_1\left(
\frac{c}{c-1}
\right)
\pmod{p},
\end{equation}
using standard congruential properties of the truncated logarithm.
Equation~\eqref{thm:2kksuk-pol*} takes place in the field $\Q_p(c)$ of rational expressions, or rational functions following the prevalent terminology;
to be more precise, it takes place in its valuation ring with respect to the $p$-adic valuation.
This is our preferred setting for congruences in this paper, as congruences between rational functions rather than polynomials.
We will extend this interpretation to algebraic functions and further clarify it after stating Theorem~\ref{thm:rkksuk_z}.

When $r=3$, finding $\mathcal{B}_r(x)$ from Equation~\eqref{eq:B_r_def}
requires solving a cubic equation.
However, in terms of the parametrization
$x=c^2/(1-c)^3$
introduced in~\cite{MatTau:truncation-cubic} one finds
$
\mathcal{B}_3(x)=(1-c)/(1-\beta),
$
where $\beta=(1-\sqrt{1-4c})/2$,
and hence Theorem~\ref{thm:rkksuk-series}
yields
\begin{equation}\label{eq:3kksuk-series}
 \sum_{k=1}^\infty\binom{3k}{k}
\frac{x^k}{k}
=3\log\left(
\frac{1-c}{1-\beta}\right)
=3\Li_1(\beta)-3\Li_1(c),
\end{equation}
an identity in the formal power series ring $\Q[[\beta]]$, where $c=\beta(1-\beta)$.
This identity has the following congruence versions, for any odd prime $p$:
\begin{equation}\label{eq:3kksuk-pol-long}
(1-c)^{2p}\sum_{0<k<p}\binom{3k}{k}\frac{x^k}{k}
\equiv
3\LL_1(\beta)-3(1-c^p)\LL_1(c)
\pmod{p},
\end{equation}
and
\begin{equation}\label{eq:3kksuk-pol-short}
(1-c)^{p}\sum_{0<k<p/3}\binom{3k}{k}\frac{x^k}{k}
\equiv
3\LL_1(\beta)-3\LL_1(c)
\pmod{p}.
\end{equation}
Here $x=c^2/(1-c)^3$ and $c=\beta(1-\beta)$,
and the congruences take place in the polynomial ring $\Z_p[\beta]$.
We originally proved Equations~\eqref{eq:3kksuk-pol-short} and~\eqref{eq:3kksuk-pol-long}
by applying the methods of~\cite{MatTau:truncation-cubic}, which were
based on truncating the power series involved and then using appropriate congruences
modulo a prime $p$.
However, those early proofs are now superseded by the the more systematic approach, for arbitrary $r$,
that we develop in Sections~\ref{sec:symmetric} and~\ref{sec:pounds}.
We will indicate in Remark~\ref{rem:beta} how to deduce Equations~\eqref{eq:3kksuk-pol-long} and~\eqref{eq:3kksuk-pol-short}
from our more general results, namely, Theorems~\ref{thm:rkksuk} and~\ref{thm:rkksuk_short}.

\section{Congruences for sums of the form $\sum_k\binom{rk}{k}x^k/k$}
\label{sec:symmetric}

In this section we present a novel approach to evaluating sums of the form $\sum_k\binom{rk}{k}x^k/k$ modulo a prime $p$,
over certain natural ranges delimited by regions where the binomial coefficients involved vanish modulo $p$.
Any prime $p$ considered in this paper will be assumed odd.

Let $q$ be a power of a prime $p$.
According to~\cite[Lemma 3.1]{MatTau:truncation-cubic},
in the range $0\le k<q$ the binomial coefficient $\binom{rk}{k}$
is a multiple of $p$ unless
$k\in A(r,m)$
for some $0<m<r$, where
\[
A(r,m)=\left\{k\in\Z:
(m-1)q/(r-1)
\le k<
mq/r
\right\}.
\]
In particular, $A(r,1)=\{k\in\Z:0\le k<q/r\}$.
In this paper we take $q=p$, and set
$A^\ast(r,m)=A(r,m)\setminus\{0\}$,
whence
$A^\ast(r,m)=A(r,m)$ unless $m=1$.

For example, while
$\sum_{0<k<p}\binom{2k}{k}x^k/k\equiv\sum_{0<k<p/2}\binom{2k}{k}x^k/k\pmod{p}$,
for the sum
$\sum_k\binom{3k}{k}x^k/k$
modulo $p$ we have two summation ranges of interest, namely,
$A(3,1)^\ast=\{k\in\Z:0<k<p/3\}$
and
$A(3,2)=\{k\in\Z:q/2\le k<2p/3\}$.
Summing over the union $A(3,1)\cup A(3,2)$ is equivalent as summing over the full range $0<k<p$.

According to Equation~\eqref{eq:B_r_equation},
the power series $z=\mathcal{B}_r(x)$ is a solution of the equation $1-z+xz^r=0$.
In the following result we consider all solutions of this equation simultaneously.
That does not depend directly on Theorem~\ref{thm:rkksuk-series} but takes inspiration from it.

\begin{theorem}\label{thm:rkksuk_z}
Let $1<r<p$, let $x,z$ be indeterminates over the $p$-adic field $\Q_p$, and let $z_i$ satisfy $1-z+xz^r=\prod_{i=1}^r(1-z/z_i)$.
Then for $1\le m<r$ we have
\[
\frac{p}{r}\sum_{k\in A^\ast(r,m)}\binom{rk}{k}\frac{x^k}{k}
\equiv
\delta(m,1)+
(-1)^m\sum_{1\leq i_1<\dots <i_m\leq r}z_{i_1}^{-p}\cdots z_{i_m}^{-p}
\pmod{p^2}.
\]
\end{theorem}

To put these and further congruences in the paper on a firm algebraic footing
while avoiding cumbersome statements, we clarify once and for all how congruences
involving indeterminates should be interpreted in this paper.

In the congruence of Theorem~\ref{thm:rkksuk_z} the symbol $x$ is to be considered an indeterminate
over the rational field $\Q$.
Therefore, $xz^r-z+1$ is a polynomial with coefficients in the field of rational functions $\Q(x)$.
Consequently, its roots $z_i$ are algebraic functions of $x$, and clearly $x=\prod_{i=1}^r(-z_i^{-1})$.
They belong to the splitting field of $xz^r-z+1$ over $\Q(x)$,
and this field extension is unramified at the prime $p$ because
we have assumed $r<p$.
This ensures that the congruence of Theorem~\ref{thm:rkksuk_z}, if desired,
can be evaluated on values of $x$ that are integers prime to $p$
(or even algebraic integers),
because then the $z_i$ and their reciprocals will evaluate to $p$-integral
algebraic numbers.
However, it will be slightly more convenient and possibly clearer to use the field $\Q_p$
of $p$-adic numbers in place of $\Q$ in our statements.

\begin{proof}
We evaluate the coefficient of $z^{mp}$ modulo $p^2$
in the polynomial $(1-z+xz^r)^p$ in two different ways, for $1\le m<r$.

On the one hand, we have
\begin{align*}
[z^{mp}](1-z+xz^r)^p
&=
\sum_{k\geq 0}\binom{p}{k}x^k[z^{mp-rk}](1-z)^{p-k}
\\&=
\sum_{k\geq 0}\binom{p}{k}x^k\binom{p-k}{mp-rk}(-1)^{mp-rk}
\\&=
\sum_{k\geq 0}\binom{p}{mp-rk}\binom{rk-(m-1)p}{k}(-1)^{mp-rk}x^k
\\&=
\sum_{k\in A(r,m)}\binom{p}{mp-rk}\binom{rk-(m-1)p}{k}(-1)^{mp-rk}x^k.
\end{align*}
Because
\[
\binom{p}{mp-rk}
=\frac{p}{mp-rk} \binom{p-1}{mp-rk-1}
\equiv \frac{(-1)^{mp-rk}}{rk}p
 \pmod{p^2}
\]
for $0<mp-rk<p$,
and because
$\binom{rk-(m-1)p}{k}\equiv\binom{rk}{k}\pmod{p}$
we conclude
\[
[z^{mp}](1-z+xz^r)^p
\equiv
-\delta(m,1)+\frac{p}{r}\sum_{k\in A^\ast(r,m)}\binom{rk}{k}\frac{x^k}{k} \pmod{p^2}.
\]

On the other hand, using the factorization of $(1-z+xz^r)^p$ we find
\begin{align*}
[z^{mp}](1-z+xz^r)^p
&=
[z^{mp}]\prod_{i=1}^r(1-z/z_i)^p\\
&=
\sum_{\substack{
0\leq j_i\leq p\\
j_1+\dots + j_r=mp
}}
\prod_{i=1}^r\binom{p}{j_i}(-z_i)^{-j_i}.
\end{align*}
In a given term of the above sum,
if at least one of $j_1,\ldots,j_r$
differs from $0$ and $p$, then at least two of them do, and hence the corresponding
binomial coefficients $\binom{p}{j_i}$ in the product are multiples of $p$.
Consequently, we find
\begin{align*}
[z^{mp}](1-z+xz^r)^p
&\equiv
\sum_{\substack{
j_i\in\{0,p\}\\
j_1+\dots + j_r=mp}}
\prod_{i=1}^r(-z_i)^{-j_i}
\pmod{p^2}
\\&\equiv
(-1)^m\sum_{1\leq i_1<\dots <i_m\leq r}z_{i_1}^{-p}\cdots z_{i_m}^{-p} \pmod{p^2}.
\end{align*}

The desired conclusion follows from comparing the result of the two calculations above.
\end{proof}

The $r-1$ congruences found in Theorem~\ref{thm:rkksuk_z} can be combined in the more compact formulation
\begin{equation*}
\frac{p}{r}\sum_{m=1}^{r-1}z^{mp}\sum_{k\in A^\ast(r,m)}\binom{rk}{k}\frac{x^k}{k}
\equiv
\prod_{i=1}^r(1-z^p/z_i^p)
-(1-z^p+x^pz^{rp})
\pmod{p^2},
\end{equation*}
where $z$ is a further indeterminate.
In particular, evaluating that for  $z=1$ we find
\begin{equation}\label{eq:rkksuk_long}
\frac{p}{r}\sum_{0<k<p}\binom{rk}{k}\frac{x^k}{k}
\equiv
-x^p+
\prod_{i=1}^r(1-z_i^{-p})
\pmod{p^2}.
\end{equation}

\section{Congruences modulo $p$ in terms of truncated logarithms}\label{sec:pounds}

In this section we apply Theorem~\ref{thm:rkksuk_z} to evaluate
the sum $\sum\binom{rk}{k}x^k/k$ modulo $p$, limiting ourselves to the ranges $0<k<p$ and $0<k<p/r$.
While Theorem~\ref{thm:rkksuk_z} is expressed in terms of the solutions of the equation
$1-z+xz^r=0$ for $z$, in the rest of the paper we find it more convenient to work with the roots of the
equation $x(c-1)^r+c^{r-1}=0$ for $c$, which also facilitates comparison with results
from~\cite{MatTau:truncation,MatTau:truncation-cubic}
when $r\le 3$.
Thus, with suitable numbering those roots are related by
$z_i=1-1/c_i$ or, equivalently,
$c_i=1/(1-z_i)$.

Several results of this paper involve instances of {\em finite polylogarithms} $\LL_s(x)=\sum_{k=1}^{p-1}x^k/k^s$,
which are truncated versions of the ordinary complex polylogarithmic series $L_s(x)=\sum_{k=1}^{\infty}x^k/k^s$
that are relevant when working modulo a prime $p$ (or modulo a power of $p$).
The parameter $s$ can be any integer, but in this paper we will only need the values $s=0,1,2$.
The finite polylogarithms satisfy various functional equations in the guise of congruences modulo $p$,
an obvious one being
\begin{equation}\label{eq:fe1}
x^p\pounds_s(1/x)\equiv(-1)^s\pounds_s(x)\pmod{p},
\end{equation}
which relates $\pounds_s(x)$ to its reciprocal polynomial and clearly has no analogue for ordinary polylogarithms.

The case of
$\pounds_0(x)=\sum_{k=1}^{p-1}x^k
=(x-x^p)/(1-x)$
will only play a marginal role in the proof of Theorem~\ref{thm:rkk}.

The more interesting {\em truncated logarithm} $\pounds_1(x)$ comes up when taking the $p$th power of a sum or difference modulo $p^2$,
because (for $p$ odd, as assumed through this paper) we have
\begin{equation}\label{eq:pth_power_mod_p^2}
(1-x)^p\equiv 1-x^p-p\pounds_1(x)\pmod{p^2}.
\end{equation}
Equation~\eqref{eq:pth_power_mod_p^2} implies the congruence
\begin{equation}\label{eq:fe2}
\pounds_1(1-x)\equiv\pounds_1(x)\pmod{p},
\end{equation}
which is another functional equation for the truncated logarithm modulo $p$.

Successive alternate applications of Equation~\eqref{eq:fe2} and the case $s=1$ of Equation~\eqref{eq:fe1}
yield six different expressions for $\pounds_1(x)$ modulo $p$,
where the argument of $\pounds_1$ is, in turn, $x, 1-x, 1/(1-x), x/(x-1), (x-1)/x, 1/x$.
In particular, one has
$\pounds_1(x)\equiv(x-1)^p\pounds_1(x/(x-1))\pmod{p}$,
which is precisely a truncated version of the functional equation
$-\log(1-x)=\log(1/(1-x))$
for the ordinary polylogarithm $L_1(x)=-\log(1-x)$.

Expanding $(1-x-y)^p$ via a suitable double application of Equation~\eqref{eq:pth_power_mod_p^2},
viewing that either as $((1-x)-y)^p$ or as $((1-y)-x)^p$ and then equating the two results,
yields a $2$-variable functional equation for the truncated logarithm, which was originally discovered by Kontsevich~\cite{Kontsevich}.
In one of its equivalent forms, that {\em $4$-term functional equation} can also be deduced from an appropriate version of the
functional equation $\log(xy)=\log(x)+\log(y)$ for the analytic logarithm, see~\cite[Equation~(31)]{MatTau:truncation}.
However, we will have no need for the $4$-term functional equation in this paper, and we refer the interested reader to~\cite{MatTau:truncation} for further discussion.

The {\em finite dilogarithm} $\pounds_2(x)$ can be used to refine Equation~\eqref{eq:pth_power_mod_p^2}
to a congruence modulo $p^3$, see~\cite[Equation~(6)]{Gra:Fermat}, namely,
\begin{equation}\label{eq:pth_power_mod_p^3}
(1-x)^p\equiv 1-x^p-p\pounds_1(x)-p^2\pounds_2(1-x)\pmod{p^3}.
\end{equation}
We will only need this congruence in the proof of Theorem~\ref{thm:rkksukk}.

Closed forms for various special values of $\pounds_1(x)$ and $\pounds_2(x)$ modulo $p$
are known, and will be recalled in Section~\ref{sec:numerical} as needed.

Now we proceed with evaluating the sum $\sum_{0<k<p} \binom{rk}{k} x^k/k$ modulo $p$ in terms of values of the truncated logarithm.
This range of summation is easier to handle than shorter ranges.

\begin{theorem}\label{thm:rkksuk}
Let $1\le r<p$, let $x,c$ be indeterminates over $\Q_p$,
and let $c_1,\ldots,c_r$ satisfy
$x(c-1)^r+c^{r-1}=x\prod_{i=1}^r(c-c_i)$.
Then
\[
\sum_{0<k<p} \binom{rk}{k} \frac{x^k}{k}
\equiv
-rx^p\sum_{i=1}^{r}
\LL_1(c_i) \pmod{p}.
\]
\end{theorem}

\begin{proof}
When $r=1$ we have $x(1-c_1)=1$.
Hence in this case the claimed congruence follows from the general congruence
$\LL_1(x)\equiv -x^p\LL_1(1/x)\pmod{p}$.
Therefore, we may assume $r>1$ in this proof, and we have $\prod_{i=1}^rc_i=1$.

According to Equation~\eqref{eq:rkksuk_long}, the claimed congruence
follows once we prove
\[
-x^p+
\prod_{i=1}^r(1-z_i^{-p})
\equiv
-px^p\sum_{i=1}^{r}
\LL_1(c_i)
\pmod{p^2},
\]
where $1/z_i=c_i/(c_i-1)$.
Indeed, we have
\begin{align*}
-x^p+
\prod_{i=1}^r(1-z_i^{-p})
&=
-x^p+x^p
\prod_{i=1}^r(1-z_i^p)
\\&\equiv
-x^p+x^p
\prod_{i=1}^r\bigl((1-z_i)^p+p\LL_1(1-z_i)\bigr)
\pmod{p^2}
\\&\equiv
-x^p+x^p
\prod_{i=1}^r\bigl((1/c_i)^p+p\LL_1(1/c_i)\bigr)
\pmod{p^2}
\\&=
-x^p+x^p
\prod_{i=1}^r\bigl(1-p\LL_1(c_i)\bigr)
\\&\equiv
-px^p\sum_{i=1}^{r}
\LL_1(c_i)
\pmod{p^2},
\end{align*}
because of Equation~\eqref{eq:pth_power_mod_p^2} and functional equations for the truncated logarithm.
\end{proof}

Dealing with the analogous sum over the range $A(r,1)$, that is, over $0<x<p/r$,
is trickier and requires a separate lemma, which we postpone as Lemma~\ref{lemma:technical} for better justification.

\begin{theorem}\label{thm:rkksuk_short}
Let $1\le r<p$, let $x,c$ be indeterminates over $\Q_p$,
and let $c_1,\ldots,c_r$ satisfy
$x(c-1)^r+c^{r-1}=x\prod_{i=1}^r(c-c_i)$.
Then
\[
\sum_{0<k<p/r}\binom{rk}{k}\frac{x^k}{k}
\equiv
-\sum_{i=1}^r\frac{\LL_1(c_i)}{(1-c_i)^p} \pmod{p}.
\]
\end{theorem}

\begin{proof}
We use the special case $m=1$ of Theorem~\ref{thm:rkksuk_z}, which reads
\begin{equation}\label{eq:rkksuk_z_short}
p\sum_{0<k<p/r}\binom{rk}{k}\frac{x^k}{k}
\equiv
r-r\sum_{i=1}^r\frac{1}{z_i^p}
\pmod{p^2},
\end{equation}
where $1/z_i=c_i/(c_i-1)$.
According to Equation~\eqref{eq:pth_power_mod_p^2}, for each $i=1,\ldots,r$ we have
\[
\frac{1}{z_i^p}
=
\left(1-\frac{1}{1-c_i}\right)^p
\equiv
1-\frac{1}{(1-c_i)^p}-p\LL_1\left(\frac{1}{1-c_i}\right)
\pmod{p^2},
\]
and
$\LL_1\bigl(1/(1-c_i)\bigr)\equiv-\LL_1(c_i)/(1-c_i)^p\pmod{p}$.
Consequently, summing over $i$ we find
\[
\sum_{i=1}^r\frac{1}{z_i^p}
\equiv r -
\sum_{i=1}^r\frac{1}{(1-c_i)^p}
+p\sum_{i=1}^r\frac{\LL_1(c_i)}{(1-c_i)^p}
\pmod{p^2}.
\]
Hence Equation~\eqref{eq:rkksuk_z_short} yields
\[
p\sum_{0<k<p/r}\binom{rk}{k}\frac{x^k}{k}
\equiv
\sum_{i=1}^r\frac{1}{(1-c_i)^p}
-(r-1)\sum_{i=1}^r\frac{1}{z_i^p}
-p\sum_{i=1}^r\frac{\LL_1(c_i)}{(1-c_i)^p}
\pmod{p^2}.
\]
Now the first two sums in the right-hand side cancel out because of Lemma~\ref{lemma:technical} below,
and the desired conclusion follows.
\end{proof}

Note that the congruence
$\sum_{i=1}^r 1/(1-c_i)^p=(r-1)\sum_{i=1}^r 1/z_i^p\pmod{p}$
follows from taking $p$th powers on each side of the identities
$\sum_{i=1}^r 1/(1-c_i)=r-1$
and $\sum_{i=1}^r 1/z_i=1$,
which in turn can be seen from
the defining equations of the quantities $1/(1-c_i)$ and $1/z_i$.
However, in order to complete the proof of Theorem~\ref{thm:rkksuk_short} we need that congruence to hold modulo $p^2$.
The fact that it does is the content of Lemma~\ref{lemma:technical}.

Our proof requires a standard fact on symmetric functions, which we recall now and we will use again in Section~\ref{sec:mod_p^2}.
If $\prod_{i=1}^r{(1+t_it)}=\sum_{j=0}^re_jt^j$,
where $e_j=e_j(t_1,\ldots,t_r)$ are the elementary symmetric functions of
the indeterminates $t_1,\ldots,t_r$, with $e_0=1$ by definition,
then the power sums $s_k=s_k(t_1,\ldots,t_n)=\sum_{i=1}^{r}t_i^k$
relate to them according to the identity
\begin{equation}\label{eq:symmetric}
\sum_{k=1}^{\infty}(-1)^{k-1}s_k\frac{t^k}{k}=\log\biggl(\sum_{j=0}^re_jt^j\biggr).
\end{equation}
In fact, this standard identity of formal power series can be found by adding together the identities
$\sum_{k=1}^{\infty}(-1)^{k-1}t_i^kt^k/k=\log(1+t_it)$
for $i=1,\ldots,r$,
and noting $\sum_{i=1}^r\log(1+t_it)=\log\bigl(\prod_{i=1}^r(1+t_it)\bigr)$.

\begin{lemma}\label{lemma:technical}
Let $1\le r<p$, let $x,c$ be indeterminates over $\Q_p$,
and let $c_1,\ldots,c_r$ satisfy
$x(c-1)^r+c^{r-1}=x\prod_{i=1}^r(c-c_i)$.
Then
\[
\sum_{i=1}^r\frac{1}{(1-c_i)^p}\equiv (r-1)\sum_{i=1}^r\frac{c_i^p}{(c_i-1)^p}\pmod{p^2}.
\]
\end{lemma}

\begin{proof}
We specialize the fact on symmetric functions recalled above taking $t_i=1/(1-c_i)$, whence
$\sum_{j=0}^{r}e_jt^j=(1+t)^{r-1}+xt^r$.
Equation~\eqref{eq:symmetric} implies
\begin{align*}
\sum_{i=1}^r\frac{1}{(1-c_i)^p}
&=p[t^p]\log\bigl((1+t)^{r-1}+xt^r\bigr)\\
&=p[t^p](r-1)\log(1+t)+p[t^p]\log\left(1+\frac{xt^r}{(1+t)^{r-1}}\right).
\end{align*}
Here we have used the standard notation $[t^k]f(t)=\Res_{t=0}\bigl(t^{-k-1}f(t)\bigr)$ for the coefficient of $t^k$ in the formal power series $f(t)$.
A similar application of Equation~\eqref{eq:symmetric} with $t_i=1/z_i$, whence
$\sum_{j=0}^{r}e_jt^j=1+t+x(-t)^r$, yields
\begin{align*}
\sum_{i=1}^r\frac{c_i^p}{(c_i-1)^p}
&=p[t^p]\log\bigl(1+t+x(-t)^r\bigr)\\
&=p[t^p]\log(1+t)+p[t^p]\log\left(1+\frac{x(-t)^r}{1+t}\right).
\end{align*}
Therefore, to complete the proof it remains to prove the congruence
\[
[t^p]\log\left(1+\frac{xt^r}{(1+t)^{r-1}}\right)
\equiv
[t^p](r-1)\log\left(1+\frac{x(-t)^r}{1+t}\right)\pmod{p}.
\]

Since $1/(1+t)^N=\sum_{j=0}^{\infty}\binom{N-1+j}{j}(-t)^j$
we have
\begin{align*}
\log\left(1+\frac{xt^r}{(1+t)^{r-1}}\right)
&=\sum_{k=1}^{\infty}\frac{(-1)^{k-1}}{k}\,\frac{(xt^r)^k}{(1+t)^{(r-1)k}}\\
&=\sum_{k=1}^{\infty}\frac{(-1)^{k-1}x^kt^{rk}}{k}\,\sum_{j=0}^{\infty}\binom{(r-1)k-1+j}{j}(-t)^j,
\end{align*}
and hence
\begin{align*}
[t^p]\log\left(1+\frac{xt^r}{(1+t)^{r-1}}\right)
&=\sum_{rk+j=p}\frac{(-1)^{k-1}x^k}{k}\binom{p-k-1}{j}(-1)^j.
\end{align*}
In a similar way we find
\begin{align*}
[t^p]\log\left(1+\frac{x(-t)^r}{1+t}\right)
(-1)^j\\
&=-\sum_{rk+j=p}\frac{(-1)^{k-1}x^k}{k}\binom{k-1+j}{j}.
\end{align*}

The desired conclusion follows because for $rk+j=p$ we have $0<j<p$, whence
\[
\binom{p-k-1}{j}(-1)^j
\equiv
\binom{-k-1}{j}(-1)^j=
\binom{k+j}{j}
\pmod{p},
\]
but also $k>0$, whence
\[
\binom{k+j}{j}=\frac{k+j}{k}\binom{k-1+j}{j}
\equiv-(r-1)\binom{k-1+j}{j}\pmod{p},
\]
where we have used $j=p-rk\equiv -rk\pmod{p}$.
\end{proof}

\begin{rem}\label{rem:mystery}
The congruence of Lemma~\ref{lemma:technical} is one of a set of three congruences modulo $p^2$
that relate either the sum of $p$th powers of all elements
$c_i$ with its analogue for the elements $1-c_i$,
or that sum for $1/(1-c_i)$ with that for $c_i/(c_i-1)$,
or that sum for $(c_i-1)/c_i$ with that for $1/c_i$.
These congruences are all rather trivial modulo $p$, but appear to be not quite obvious modulo $p^2$,
and can be proved to hold in a similar way as the proof of Lemma~\ref{lemma:technical}.
Thus, one can show
\begin{equation}\label{eq:mystery2}
(r-1)\sum_{i=1}^{r} c_i^p+r\sum_{i=1}^{r} (1-c_i)^p
\equiv x^{-p}+r(r-1)
\pmod{p^2},
\end{equation}
whose version modulo $p$ follows from taking $p$th powers of each side of the identities
$\sum_{i=1}^rc_i=r-1/x$,
$\sum_{i=1}^r(1-c_i)=1/x$.
Similarly, for $r>2$ one can show
\[
\sum_{i=1}^{r} \frac{1}{c_i^p}+r\sum_{i=1}^{r} \left(1-\frac{1}{c_i}\right)^p
\equiv
r
\pmod{p^2};
\]
when $r=2$ the congruence is similar but the right-hand side is $x^{-p}+2$ instead.
\end{rem}

By differentiating the congruences of Theorems~\ref{thm:rkksuk} and~\ref{thm:rkksuk_short} we infer corresponding congruences for sums $\sum_{k} \binom{rk}{k} x^k$.

\begin{theorem}\label{thm:rkk}
Let $1\le r<p$, let $x,c$ be indeterminates over $\Q_p$,
and let $c_1,\ldots,c_r$ satisfy
$x(c-1)^r+c^{r-1}=x\prod_{i=1}^r(c-c_i)$.
Then
\[
\sum_{0<k<p} \binom{rk}{k} x^k
\equiv
-rx^p\sum_{i=1}^{r}
\frac{c_i-c_i^p}{r-1+c_i} \pmod{p}
\]
and
\[
\sum_{0<k<p/r} \binom{rk}{k} x^k
\equiv
-\sum_{i=1}^r\frac{c_i-c_i^p}{(1-c_i^p)(r-1+c_i)} \pmod{p}.
\]
\end{theorem}

\begin{proof}
We obtain the desired conclusions by differentiating the congruences of
Theorems~\ref{thm:rkksuk} and~\ref{thm:rkksuk_short} with respect to $x$ and then multiplying by $x$.
Because $x=-c_i^{r-1}/(c_i-1)^r$ for each $i=1,\ldots,r$, we have
\[
\frac{1}{x}\frac{dx}{dc_i}=\frac{r-1}{c_i}-\frac{r}{c_i-1}
=\frac{r-1+c_i}{c_i(1-c_i)}.
\]
Note that such derivatives can be taken in a formal way, as we are dealing with algebraic functions,
but this makes no difference from taking them in the $p$-adic analytic sense.
Since $d\pounds_s(c_i)/dc_i=d\pounds_{s-1}(c_i)/c_i$ for every integer $s$,
where we have conveniently set
$\pounds_0(x)=\sum_{k=1}^{p-1}x^k=(x-x^p)/(1-x)$,
we find
\[
x\frac{d\pounds_s(c_i)}{dx}
=\frac{x}{dx/dc_i}\frac{d\pounds_s(c_i)}{dc_i}
=\frac{1-c_i}{r-1+c_i}\pounds_{s-1}(c_i).
\]
Using this with $s=1$, Theorems~\ref{thm:rkksuk} and~\ref{thm:rkksuk_short} imply
\[
\sum_{0<k<p} \binom{rk}{k} x^k
\equiv
-rx^p\sum_{i=1}^{r}
\frac{1-c_i}{r-1+c_i}
\LL_0(c_i) \pmod{p}
\]
and
\[
\sum_{0<k<p/r} \binom{rk}{k} x^k
\equiv
-\sum_{i=1}^r\frac{1-c_i}{(1-c_i)^p(r-1+c_i)}\LL_0(c_i) \pmod{p}.
\]
The desired conclusions follow since
$\pounds_0(c_i)=\sum_{k=1}^{p-1}c_i^k
=(c_i-c_i^p)/(1-c_i)$.
\end{proof}

If we wish to substitute an algebraic integer $a$ for the indeterminate $x$
in the congruences of Theorem~\ref{thm:rkk},
or, more conveniently, an element $a$ of some finite extension of $\Q_p$
with $v_p(a)\ge 0$,
a couple of degenerate cases need to be avoided.
Firstly, we should assume $a\not\equiv 0\pmod{p}$,
otherwise the sums to be evaluated would be zero modulo $p$ for a trivial reason.
Secondly, the denominators $r-1+c_i$ that appear in the right-hand sides
of those congruences should also be nonzero modulo $p$, and this means
$a\not\equiv (r-1)^{r-1}/r^r\pmod{p}$.
Because the discriminant of the polynomial $a(c-1)^r+c^{r-1}$ equals
$(-1)^{r(r-1)/2}a^{r-2}(r^ra-(r-1)^{r-1})$,
such restrictions on $a$ also ensure that the polynomial viewed modulo $p$
has $r$ distinct roots in a splitting field over the field $\F_p$ of $p$ elements.
If all those roots $c_1,\ldots,c_r$ actually belong to $\F_p$, then
because $c_i^p=c_i$ we obtain the following interesting consequence.

\begin{cor}\label{cor:rkk}
Let $1\le r<p$, and let $a\in\Z_p$
be such that
$a\not\equiv 0,(r-1)^{r-1}/r^r\pmod{p}$.
If the polynomial $a(c-1)^r+c^{r-1}\in\F_p[c]$ factors completely in $\F_p[c]$, then
$\sum_{0<k<p}\binom{rk}{k}a^k\equiv 0\pmod{p}$
and
$\sum_{0<k<p/r}\binom{rk}{k}a^k\equiv 0\pmod{p}$.
\end{cor}

Because the congruences of Theorem~\ref{thm:rkk} were already known in essence if $r\le 3$, a number of comments and some questions arise.

\begin{rem}
For $r\le 3$ the congruences of Theorem~\ref{thm:rkk} were known from~\cite{MatTau:truncation} and~\cite{MatTau:truncation-cubic},
but see also~\cite{SunZH:sums_(3k_k)} and~\cite{Pilehrood15},
and in a more general form that allowed for a power $q$ of $p$ in place of $p$ in the summation ranges.
In particular, when formulated in the setting of Theorem~\ref{thm:rkk},
the congruences of~\cite[Corollary~3.4]{MatTau:truncation-cubic}
become
\begin{equation}\label{eq:3kk-pol-long}
\sum_{0<k<q}\binom{3k}{k}x^k
\equiv
-3x^q\sum_{i=1}^3\frac{c_i-c_i^q}{2+c_i}
\pmod{p},
\end{equation}
\begin{equation}\label{eq:3kk-pol-short}
\sum_{0<k<q/3}\binom{3k}{k}x^k
\equiv
-\sum_{i=1}^3\frac{c_i-c_i^q}{(1-c_i^q)(2+c_i)}
\pmod{p}.
\end{equation}
This prompts a natural question as to whether the congruences of Theorem~\ref{thm:rkk} also hold with $p$ appropriately replaced by $q$.
It seems that such an extension cannot be attained by the approach used here,
which depends on Theorems~\ref{thm:rkksuk} and~\ref{thm:rkksuk_short}, where replacing $p$ with $q$ does not make sense.
\end{rem}

\begin{rem}\label{rem:MT25}
The congruences of~\cite[Corollary~3.4]{MatTau:truncation-cubic} were actually expressed in terms of a single indeterminate $c_1$
(written as $c$ there),
rather than the three algebraic functions $c_1,c_2,c_3$ of $x$.
The reason why that worked is that after writing $c_1$ as $c_1=\beta(1-\beta)$, we may express $c_2$ and $c_3$ as
$c_2=-c_1/\beta^3$ and $c_3=-c_1/(1-\beta)^3$.
If we do that in the congruences of Theorem~\ref{thm:rkk} and work out the sums,
we can arrange the results in a way that the only appearance of $\beta$ occurs in an expression $(1-2\beta)^{p-1}$,
which equals $(1-4c_1)^{(p-1)/2}$, so that $\beta$ can be eliminated altogether from the congruences.
\end{rem}

\begin{rem}\label{rem:beta}
In a similar fashion, when $r=3$ the congruences of Theorems~\ref{thm:rkksuk} and~\ref{thm:rkksuk_short}
can be written in terms of $c_1$ and $\beta$ only, and become our Equations~\eqref{eq:3kksuk-pol-long} and~\eqref{eq:3kksuk-pol-short}.
We have presented those in Section~\ref{sec:3ksuk_series} (using the letter $c$ in place of $c_1$)
because of their striking similarity with the generating function of the corresponding series.
Note that we have not explicitly eliminated  $\beta$ from Equations~\eqref{eq:3kksuk-pol-long} and~\eqref{eq:3kksuk-pol-short}.
However, $\beta$ only occurs in the expression $\LL_1(\beta)$,
which is congruent modulo $p$ to a certain polynomial in $c_1=\beta(1-\beta)$,
because $\LL_1(1-\beta)\equiv\LL_1(\beta)\pmod{p}$.
Thus, Equations~\eqref{eq:3kksuk-pol-long} and~\eqref{eq:3kksuk-pol-short} may be viewed as congruences between polynomials
in the indeterminate $c_1$ (written as $c$ there) rather than in $\beta$.

Alternately, Equations~\eqref{eq:3kksuk-pol-long} and~\eqref{eq:3kksuk-pol-short}
can also be proved by integrating
Equations~\eqref{eq:3kk-pol-long} and~\eqref{eq:3kk-pol-short}
with $q=p$,
after those are written in terms of polylogarithms
$\pounds_0(\beta)$ and $\LL_0(c)$.
However, because those finite polylogarithms are actually rational functions,
this rewriting was not easy to recognize (and was inessential) in the context
of~\cite{MatTau:truncation-cubic}.
\end{rem}

\section{Congruences modulo $p^2$}\label{sec:mod_p^2}

In this section we introduce a different method that allows us
to strengthen the congruence modulo $p$ of Theorem~\ref{thm:rkksuk} to a congruence modulo $p^2$,
which appears in Theorem~\ref{thm:rkksukmod2}.\
However, this method only appears to works smoothly for  sums over the whole range $0<k<p$,
and not for the shorter range $0<k<p/r$.
The method relies on certain identities that we state in Theorem~\ref{thm:identity}.

We will need the standard fact on symmetric functions that we recalled before Lemma~\ref{lemma:technical}
and is summarized by Equation~\eqref{eq:symmetric}.
We specialize the indeterminates $t_i$
to the elements $c_i$ introduced in the previous section,
which relate to $x$ according to
$x(c-1)^r+c^{r-1}=x\prod_{i=1}^r(c-c_i)$.
However, our identities will be best stated and proved in terms of
an auxiliary indeterminate $y=1/x$.

\begin{theorem}\label{thm:identity}
Let $y,t$ be indeterminates over $\Q$,
and let $c_1,\ldots,c_r$ be the algebraic functions of $y$ determined by
$\prod_{i=1}^r(1+c_it)=(1+t)^r-yt$.
Set $s_k=\sum_{i=1}^r c_i^k$ and let $s_k'$ denote the derivative $ds_k/dy$.
Then we have the identities
\begin{equation}\label{id0}
\sum_{k=0}^{n-1}
\binom{rk}{k}y^{n-k-1}=\sum_{k=1}^n\binom{rn}{n-k}(-1)^k\frac{s_k'}{k},
\end{equation}
\begin{equation}\label{id1b}
\sum_{k=0}^{n-1}
\binom{rk}{k}\frac{y^{n-k}}{n-k}
=\sum_{k=1}^n\binom{rn}{n-k}(-1)^k\frac{s_k-r}{k},
\end{equation}
\begin{equation}\label{id2b}
\begin{split}
\sum_{k=0}^{n-1}
\binom{rk}{k}\frac{y^{n-k}}{(n-k)^2}
&=
-(r-1)\sum_{k=1}^n\binom{rn}{n-k}(-1)^k\frac{s_k-r}{k^2}\\
&\quad
+r\sum_{k=1}^n\binom{rn}{n-k}\frac{(-1)^k}{k}
\sum_{j=1}^k\frac{s_j-r}{j}.
\end{split}
\end{equation}
\end{theorem}

In the special case $r=2$ the three identities of Theorem~\ref{thm:identity}
are related to~\cite[Equations~(34),~(35),~(36)]{MatTau:polylog}.
However, those identities were expressed in terms of certain polynomial
Lucas sequences in place of
the sums $s_k$ and their derivatives used here.

\begin{proof}
We start with proving Equation~\eqref{id0}, after which the other two identities will follow by suitable integration.
According to Equation~\eqref{eq:symmetric} we have
\[
\sum_{k=1}^{\infty}(-1)^k\frac{s_k}{k}t^k
=-\log\bigl((1+t)^r-yt\bigr).
\]
Taking derivatives with respect to $y$ we find
\begin{align*}
\sum_{k=1}^{\infty}(-1)^k\frac{s_k'}{k}t^k
&=
-\frac{d}{dy}\log\bigl((1+t)^r-yt\bigr)
=
\frac{t}{(1+t)^r-yt}
\\&=
\frac{t}{(1+t)^r}\cdot \left(1-\frac{yt}{(1+t)^r}\right)^{-1}
=
\sum_{j=1}^{\infty}\frac{t^j}{(1+t)^{rj}}\, y^{j-1}
\\&=
\sum_{j=1}^{\infty}\sum_{i=0}^{\infty}\binom{-rj}{i}t^{i+j}\, y^{j-1}
=
\sum_{k=1}^{\infty}t^k\sum_{j=1}^{k}\binom{-rj}{k-j}\, y^{j-1}.
\end{align*}
Consequently, for $k\ge 1$ we have
\begin{equation}\label{eq:s_k'}
(-1)^k\frac{s_k'}{k}
=\sum_{j=1}^{k}\binom{-rj}{k-j} y^{j-1}.
\end{equation}
Substituting this into the right-hand side of Equation~\eqref{id0} we find
\begin{align*}
\sum_{k=1}^{n}\binom{rn}{n-k}(-1)^k\frac{s_k'}{k}
&=
\sum_{k=1}^{n}\binom{rn}{n-k}\sum_{j=1}^{k}\binom{-rj}{k-j}\, y^{j-1}
\\&=
\sum_{j=1}^{n}y^{j-1}\sum_{k=j}^{n}\binom{rn}{n-k}\binom{-rj}{k-j}
\\&=
\sum_{j=1}^{n}y^{j-1}\binom{r(n-j)}{n-j}=\sum_{k=0}^{n-1}\binom{rk}{k}y^{n-1-k},
\end{align*}
which is the left-hand side of Equation~\eqref{id0}.

To prove Equation~\eqref{id1b}, note that differentiating each side with respect to $y$ returns Equation~\eqref{id0}.
Also, Equation~\eqref{id1b} holds when setting $y=0$ because then $c_1=\dots=c_r=1$,
and hence $s_k=r$.

Similarly, to prove Equation~\eqref{id2b} it suffices to show that differentiating each side and then multiplying by $y$ returns Equation~\eqref{id1b}.
In fact, doing that to its right-hand side (after replacing the summation index $j$ with $i$) and then using Equation~\eqref{eq:s_k'} yields
\[
\sum_{k=1}^n\binom{rn}{n-k}
\sum_{j=1}^k\left(
\frac{1-r}{k}
\binom{-rj}{k-j}
+\frac{(-1)^kr}{k}
\sum_{i=j}^k(-1)^i\binom{-rj}{i-j}
\right)y^j.
\]
Now the quantity in parentheses equals
\begin{align*}
\frac{1-r}{k}
\binom{-rj}{k-j}
+\frac{r}{k}
\binom{-rj-1}{k-j}
&=
\left(\frac{1-r}{k}
+\frac{(r-1)j+k}{kj}\right)
\binom{-rj}{k-j}\\
&=
\frac{1}{j}
\binom{-rj}{k-j}.
\end{align*}
Since integrating Equation~\eqref{eq:s_k'} yields
\[
(-1)^k\frac{s_k-r}{k}
=\sum_{j=1}^{k}\binom{-rj}{k-j}\frac{y^j}{j},
\]
the desired conclusion follows.
\end{proof}

We are going to set $n=p$ in each of the identities of Theorem~\ref{thm:identity}
in order to deduce certain congruences.
Although our main goal is an evaluation modulo $p^2$ of the sum
$\sum_{k=1}^{p-1}\binom{rk}{k}x^k/k$
in Theorem~\eqref{thm:rkksukmod2},
that cannot be deduced from Equation~\eqref{id1b} alone because of the need to
evaluate $1/(p-k)$ modulo $p^2$,
but requires a contribution from
Equation~\eqref{id2b}.
Thus, we start with deducing an appropriate congruence from Equation~\eqref{id2b}.
We need $p>3$ because we will use the congruence $\LL_1(1)\equiv 0\pmod{p^2}$
in the proof of Theorem~\ref{thm:rkksukk},
and also $\LL_2(1)\equiv 0\pmod{p}$
in the proof of Theorem~\ref{thm:rkksukmod2}.

\begin{theorem}\label{thm:rkksukk}
Let $p>3$ be a prime and let $r\geq 1$.
Let $x,c$ be indeterminates over $\Q_p$,
and let $c_1,\ldots,c_r$ satisfy
$x(c-1)^r+c^{r-1}=x\prod_{i=1}^r(c-c_i)$.
Then
\begin{align*}
\sum_{0<k<p}
\binom{rk}{k}\frac{x^k}{k^2}
&\equiv
\frac{1}{p^2}
\biggl(-1+(r-1)x^p\sum_{i=1}^r(c_i^p-1)+rx^p\sum_{i=1}^r(1-c_i)^p\biggr)
\\&\quad
+rx^p\sum_{i=1}^r\LL_2(1-c_i)
\pmod{p}.
\end{align*}
\end{theorem}

\begin{proof}
Setting $n=p$ and $y=1/x$ in Equation~\eqref{id2b}, and conveniently separating the terms where multiples of $p^2$ appear as denominators, we find
\begin{align*}
\frac{x^{-p}}{p^2}+\sum_{k=1}^{p-1}
\binom{rk}{k}\frac{x^{k-p}}{(p-k)^2}
&=
(r-1)\frac{s_p-r}{p^2}
-(r-1)\sum_{k=1}^{p-1}\binom{rp}{p-k}(-1)^k\frac{s_k-r}{k^2}\nonumber\\
&\qquad
-r\frac{s_p-r}{p^2}-\frac{r}{p}\sum_{j=1}^{p-1}\frac{s_j-r}{j}
\\&\qquad
+r\sum_{k=1}^{p-1}\binom{rp}{p-k}\frac{(-1)^k}{k}
\sum_{j=1}^k\frac{s_j-r}{j}.
\end{align*}
Now multiply this identity by $p^2x^p$ and view the result modulo $p^3$.
Because
$\binom{rp}{p-k}$ is a multiple of $p$ in the range considered, because
\[
\sum_{j=1}^{p-1}\frac{s_j-r}{j}
=\sum_{i=1}^r\bigl(\LL_1(c_i)-\LL_1(1)\bigr),
\]
and because $\LL_1(1)=\sum_{j=1}^{p-1}1/j\equiv 0\pmod{p^2}$ if $p>3$ according to Wolstenholme's theorem,
we are left with
\[
p^2\sum_{k=1}^{p-1}
\binom{rk}{k}\frac{x^k}{k^2}
\equiv (-1-s_px^p+rx^p)
-prx^p\sum_{i=1}^r\LL_1(c_i)\pmod{p^3}.
\]
Using Equation~\eqref{eq:pth_power_mod_p^3} we find
\[
-p\sum_{i=1}^r\LL_1(c_i)
\equiv
\biggl(s_p+\sum_{i=1}^r(1-c_i)^p-r\biggr)
+p^2\sum_{i=1}^r\LL_2(1-c_i)\pmod{p^3},
\]
and the desired congruence follows.
\end{proof}

In the next result we use Equation~\eqref{id2b} and Theorem~\ref{thm:rkksukk} to prove a congruence modulo $p^2$.

\begin{theorem}\label{thm:rkksukmod2}
Let $p>3$ be a prime and let $r\geq 1$.
Let $x,c$ be indeterminates over $\Q_p$,
and let $c_1,\ldots,c_r$ satisfy
$x(c-1)^r+c^{r-1}=x\prod_{i=1}^r(c-c_i)$.
Then
\begin{align*}
\sum_{0<k<p}\binom{rk}{k}\frac{x^k}{k}
&\equiv
\frac{1}{p}
\biggl(2-(r-2)x^p\sum_{i=1}^r(c_i^p-1)-rx^p\sum_{i=1}^r(1-c_i)^p\biggr)
\\&\quad
+prx^p\sum_{i=1}^r\bigl(\LL_2(c_i)-\LL_2(1-c_i)\bigr)
\pmod{p^2}.
\end{align*}
\end{theorem}

\begin{proof}
Setting $n=p$ and $y=1/x$ in Equation~\eqref{id1b} we find
\[
\frac{1}{p}+\sum_{k=1}^{p-1}\binom{rk}{k}\frac{x^k}{p-k}
=-\frac{x^p(s_p-r)}{p}+x^p\sum_{k=1}^{p-1}\binom{rp}{p-k}(-1)^k\frac{s_k-r}{k}.
\]
Because
$
1/(k-p)\equiv 1/k+p/k^2\pmod{p^2}
$,
and
\[
(-1)^k\binom{rp}{p-k}
\equiv-\frac{rp}{k}\pmod{p^2}
\]
for $0<k<p$, we deduce the congruence
\[
\sum_{k=1}^{p-1}
\binom{rk}{k}\frac{x^k}{k}+p\sum_{k=1}^{p-1}
\binom{rk}{k}\frac{x^k}{k^2}\equiv \frac{1+x^p(s_p-r)}{p}+prx^p\sum_{k=1}^{p-1}
\frac{s_k-r}{k^2}
\pmod{p^2}.
\]
Because $\pounds_2(1)=\sum_{k=1}^{p-1}1/k^2\equiv 0\pmod{p}$ if $p>3$, the congruence amounts to
\begin{align*}
\sum_{k=1}^{p-1}
\binom{rk}{k}\frac{x^k}{k}
&\equiv
\frac{1}{p}\biggl(1+x^p\sum_{i=1}^{r}(c_i^p-1)\biggr)
\\&\quad
-p\sum_{k=1}^{p-1}\binom{rk}{k}\frac{x^k}{k^2}
+rpx^p\sum_{i=1}^{r}\pounds_2(c_i)\pmod{p^2},
\end{align*}
and then the conclusion follows by applying Theorem~\ref{thm:rkksukk}.
\end{proof}

Next, we deduce from Equation~\eqref{id0} a version modulo $p^2$ of the first congruence of
Theorem~\ref{thm:rkk}.

\begin{theorem}\label{thm:rkkmod2}
Let $p>2$ be a prime and let $r\geq 1$.
Let $x,c$ be indeterminates over $\Q_p$,
and let $c_1,\ldots,c_r$ satisfy
$x(c-1)^r+c^{r-1}=x\prod_{i=1}^r(c-c_i)$.
Then
\[
\sum_{0\le k<p}\binom{rk}{k}x^k
\equiv
-x^{p}\sum_{i=1}^r\frac{c_i-1}{r-1+c_i}\bigl(r-(r-1)c_i^p -r(1-c_i)^p\bigr)
\pmod{p^2}.
\]
\end{theorem}

\begin{proof}
Starting from Equation~\eqref{id0} and decomposing its right-hand side into partial fractions we find
\begin{align*}
\sum_{k=0}^{n-1}\binom{rk}{k}y^{n-k-1}&=\sum_{k=1}^n\binom{rn}{n-k}(-1)^k\frac{s_k'}{k}\\
&=\sum_{k=1}^n\binom{rn}{n-k}[t^k]\frac{t}{(1+t)^r-yt}\\
&=\sum_{k=1}^n\binom{rn}{n-k}[t^k]\sum_{i=1}^r\frac{A_i}{1+t_it}\\
&=\sum_{i=1}^rA_i\sum_{k=1}^n\binom{rn}{n-k}(-c_i)^k
\end{align*}
where $(1+t)^r-yt=\prod_{i=1}^r(1+c_it)$ and
\[
A_i=c_i\,\Res_{t=-1/c_i}\left(\frac{t}{(1+t)^r-yt}\right)
=\frac{c_i-1}{y(r-1+c_i)}.
\]
On setting $n=p$ and $y=1/x$ we find
\begin{align*}
\sum_{k=0}^{p-1}
\binom{rk}{k}x^k&=x^{p}\sum_{i=1}^r\frac{c_i-1}{r-1+c_i}\sum_{k=1}^p\binom{rp}{p-k}(-c_i)^k\\
&\equiv x^{p}\sum_{i=1}^r\frac{c_i-1}{r-1+c_i}\biggl(-c_i^p-pr\sum_{k=1}^{p-1}\frac{c_i^k}{k}\biggr)
\pmod{p^2}
\\&\equiv
-x^{p}\sum_{i=1}^r\frac{c_i-1}{r-1+c_i}\bigl(c_i^p +pr\LL_1(c_i)\bigr) \pmod{p^2},
\end{align*}
which is equivalent to the desired conclusion according to Equation~\eqref{eq:pth_power_mod_p^2}.
\end{proof}

Another route to a congruence similar to that of Theorem~\ref{thm:rkkmod2} is
applying the differential operator $x\,d/dx$ to both sides of the congruence of Theorem~\ref{thm:rkksukmod2}.
However, this procedure does not yield quite the same congruence of Theorem~\ref{thm:rkkmod2}, but rather the following variant.
The equivalence of the congruence of Theorem~\ref{thm:rkkmod2_var} with the congruence of Theorem~\ref{thm:rkkmod2} is not quite obvious,
as it depends on Equation~\eqref{eq:mystery2}.

\begin{theorem}\label{thm:rkkmod2_var}
Let $p>2$ be a prime and let $r\geq 1$.
Let $x,c$ be indeterminates over $\Q_p$,
and let $c_1,\ldots,c_r$ satisfy
$x(c-1)^r+c^{r-1}=x\prod_{i=1}^r(c-c_i)$.
Then
\[
\sum_{0<k<p}\binom{rk}{k}x^k
\equiv
x^p\sum_{i=1}^r\frac{r}{r-1+c_i}\bigl(r-1-(r-1)c_i^p-r(1-c_i)^p\bigr)\pmod{p^2}.
\]
\end{theorem}

Note that the summation range excludes $k=0$ here.
\begin{proof}
Since
$\ds
x\,dc_i/dx
=c_i(1-c_i)/(r-1+c_i)
$
as in the Proof of Theorem~\ref{thm:rkk}, we have
\[
x\,\frac{d}{dx}c_i^p
=p\frac{1-c_i}{r-1+c_i}\cdot c_i^p,
\quad
x\,\frac{d}{dx}(1-c_i)^p
=-p\frac{c_i}{r-1+c_i}\cdot (1-c_i)^p,
\]
and
\[
x\,\frac{d}{dx}\LL_2(c_i)
=\frac{1-c_i}{r-1+c_i}\cdot \LL_1(c_i),
\quad
x\,\frac{d}{dx}\LL_2(1-c_i)
=-\frac{c_i}{r-1+c_i}\cdot \LL_1(1-c_i).
\]
By applying the differential operator $x\,d/dx$ to both sides of the congruence of Theorem~\ref{thm:rkksukmod2}
and using $\LL_1(c_i)\equiv\LL_1(1-c_i)\pmod{p}$ we obtain
\begin{align*}
\sum_{k=1}^{p-1}\binom{rk}{k}x^k
&\equiv
-(r-2)x^p\sum_{i=1}^r(c_i^p-1)-rx^p\sum_{i=1}^r(1-c_i)^p
\\&\quad
-(r-2)x^p\sum_{i=1}^r\frac{1-c_i}{r-1+c_i}c_i^p
+rx^p\sum_{i=1}^r\frac{c_i}{r-1+c_i}(1-c_i)^p
\\&\quad
+px^p\sum_{i=1}^r\frac{r}{r-1+c_i}\LL_1(c_i)
\\&\equiv
r(r-2)x^p
-rx^p\sum_{i=1}^r\frac{r-2}{r-1+c_i}c_i^p
-rx^p\sum_{i=1}^r\frac{r-1}{r-1+c_i}(1-c_i)^p
\\&\quad
+px^p\sum_{i=1}^r\frac{r}{r-1+c_i}\LL_1(c_i).
\end{align*}
The desired congruence follows since
$\sum_{i=1}^r 1/(r-1+c_i)=1$.
\end{proof}

The presence of denominators $r-1+c_i$ in the congruences of Theorems~\ref{thm:rkkmod2} and  Theorems~\ref{thm:rkkmod2_var}
prevents a specialization of the congruence where any of the numbers $c_i$ takes the value $1-r$.
That occurs precisely when $x$ takes the value $x_0=(r-1)^{r-1}/r^r$,
and is equivalent to the vanishing of the discriminant of $x(c-1)^r+c^{r-1}$ as a polynomial in $c$.
As one can see using the derivative criterion, in that case the polynomial has $1-r$ as a double root, and the remaining roots are simple.

Because the congruences obtained in~\cite[Section~5]{MatTau:truncation-cubic} for the special case $r=3$
depended on just one of the elements $c_i$, this offered a workaround to evaluating on $x_0=4/27$ by selecting the simple root,
and so in~\cite[Section~5]{MatTau:truncation-cubic} we found
\[
\sum_{0\le k<q} \binom{3k}{k} (4/27)^k
\equiv 1/9\pmod{p},
\quad
\sum_{0\le k<q/3} \binom{3k}{k} (4/27)^k
\equiv 1/3\pmod{p},
\]
the former of which was already known from~\cite[Theorem 3.1]{Sun:sums_higher_Catalan}.
Since that workaround is not available in the present setting, we cover the special case $x_0=(r-1)^{r-1}/r^r$ with the following result.

\begin{theorem}\label{thm:rkkmod2_multiple}
Let $p>3$ be a prime, let $r>1$. Set $x_0=(r-1)^{r-1}/r^r$ and let
$x_0(c-1)^r+c^{r-1}=x_0\prod_{i=1}^r(c-c_i)$,
with $c_1=c_2=1-r$.
Then
\begin{align*}
\sum_{0\le k<p}\binom{rk}{k}x_0^k
&\equiv
\frac{2x_0^p}{3}\left((r-2+3pr)(r-1)^{p-1} -\frac{pr(r-2)}{r-1}\LL_1(1-r)\right)\\
&\quad -x_0^{p}\sum_{i=3}^r\frac{c_i-1}{r-1+c_i}\left(c_i^p +rp\LL_1(c_i)\right)
\pmod{p^2}.
\end{align*}
\end{theorem}

\begin{proof}
As discussed above, if $x_0=(r-1)^{r-1}/r^r$ then $c_1=c_2=1-r$ and the other numbers $c_i$ are distinct, and different from $1-r$.
The proof is similar to that of Theorem~\ref{thm:rkkmod2},
except that the partial fraction decomposition for the double root $c_1=c_2$
differs from that for the remaining roots.
One finds
\begin{align*}
\sum_{k=0}^{n-1}\binom{rk}{k}y^{n-k-1}
&=
B\sum_{k=1}^n\binom{rn}{n-k}(r-1)^k
+C\sum_{k=1}^n\binom{rn}{n-k}(k+1)(r-1)^k
\\&\quad
+\sum_{i=3}^rA_i\sum_{k=1}^n\binom{rn}{n-k}(-c_i)^k,
\end{align*}
where
\[
B=-\frac{4(r+1)}{3y(r-1)},
\qquad
C=\frac{2r}{y(r-1)},
\qquad
A_i
=\frac{c_i-1}{y(r-1+c_i)}.
\]
On setting $n=p$ we find
\begin{align*}
\sum_{k=1}^p\binom{rp}{p-k}(r-1)^k
&\equiv
(r-1)^p-rp\sum_{k=1}^{p-1}\frac{(1-r)^k}{k}
\pmod{p^2}
\\&=
(r-1)^p
-rp\LL_1(1-r),
\end{align*}
and
\begin{align*}
\sum_{k=1}^p\binom{rp}{p-k}k(r-1)^k
&\equiv
-p(1-r)^p-rp\sum_{k=1}^{p-1}(1-r)^k
\pmod{p^2}.
\end{align*}
Putting all pieces together we find the desired conclusion.
\end{proof}

\section{Numerical applications}\label{sec:numerical}

In this final section we collect various applications to numerical congruences, which we obtain by evaluating some of the generic congruences
of Sections~\ref{sec:pounds} and~\ref{sec:mod_p^2} to various values of $x$.
Several of those will involve values of the {\em Fermat quotient} $q_p(x)=(x^{p-1}-1)/p$.
Although this object is traditionally considered for $x$ an integer prime to $p$, in which case $q_p(x)$ is an integer as well,
if we regard it as a polynomial in the indeterminate $x$ it relates nicely to the truncated logarithm via the congruence
    \begin{equation}\label{eq:pth_power_mod_p^2_var}
\LL_1(x)\equiv\frac{1-x^p-(1-x)^p}{p}=-xq_p(x)-(1-x)q_p(1-x)\pmod{p},
\end{equation}
which is a rewriting of Equation~\eqref{eq:pth_power_mod_p^2}.

In our first application of our results we take $r=3$ in Theorems~\ref{thm:rkksukk},~\ref{thm:rkksukmod2},~\ref{thm:rkkmod2} and evaluate those congruences at $x=2$.
Then up to order we have $c_1=1/2$, $c_2=1+i$, $c_3=1-i$,
and so we need to evaluate the various quantities involved in the congruences for such values.
Thus, we have
\[
x^p\sum_{i=1}^r(1-c_i)^p
=2^p\bigl(2^{-p}+(-i)^p+i^p\bigr) =1.
\]
Because $(\pm i)^p=\pm(-1)^{\frac{p-1}{2}}i$ if $p$ is odd, for every positive integer $n$ we find
\begin{align*}
(1\pm i)^p
&=\leg{2}{p} 2^{\frac{p-1}{2}}\bigl(1\pm (-1)^{\frac{p-1}{2}}i\bigr)
\\&\equiv
\bigl(1\pm (-1)^{\frac{p-1}{2}}i\bigr)\sum_{k=0}^{n-1}\binom{1/2}{k}(p\,q_p(2))^k\pmod{p^n},
\end{align*}
where we have used~\cite[Lemma~4.1]{MatTau:polylog} for the congruence.
By applying this with $n=3$ we find
\begin{align*}
x^p\sum_{i=1}^r(c_i^p-1)
&=2^p\bigl(2^{-p}+(1+i)^p+(1-i)^p-3\bigr)\\
&\equiv 1+2\bigl(1+p\,q_p(2)\bigr)\biggl(2\sum_{k=0}^{2}\binom{1/2}{k}(p\,q_p(2))^k-3\biggr)
\pmod{p^3}\\
&\equiv -1+\frac{3}{2}p^2\,q_p(2)^2\pmod{p^3}.
\end{align*}
According to~\cite[Section~4]{MatTau:polylog} the values of $\pounds_2$ involved in our congruences can be evaluated modulo $p$ as
\begin{align*}
\pounds_2(1/2)&\equiv-\frac{1}{2}\,q_p(2)^2\pmod{p},
\\
\pounds_2(1\pm i)&\equiv-\frac{1}{8}\,q_p(2)^2\bigl(1\pm(-1)^{\frac{p-1}{2}}i\bigr)+\frac{1}{2}(-1)^{\frac{p-1}{2}}\,E_{p-3}\pmod{p},
\\
\pounds_2(\pm i)&\equiv \frac{1}{2}\bigl((-1)^{\frac{p-1}{2}}\pm i\bigr)\,E_{p-3}\pmod{p},
\end{align*}
where $E_{p-3}$ denotes an Euler number.
Now Theorem \ref{thm:rkksukmod2} yields the congruence
\[
\sum_{0<k<p}
\binom{3k}{k}\frac{2^k}{k}\equiv \frac{3}{2}p\,q_p(2)^2-\frac{9}{2}p\,q_p(2)^2
=-3p\,q_p(2)^2\pmod{p^2}.
\]
This is~\cite[Conjecture~9]{ConjSZW},
and a version modulo $p$ was proved in~\cite{Zhao-Pan-Sun}.
Theorem \ref{thm:rkksukk} yields
\[
\sum_{0<k<p}
\binom{3k}{k}\frac{2^k}{k^2}\equiv
6(-1)^{\frac{p-1}{2}}\,E_{p-3}\pmod{p},
\]
and for $p\not=5$ Theorem \ref{thm:rkkmod2} yields
\[
\sum_{0\le k<p} \binom{3k}{k}2^k
\equiv \frac{6(-1)^{\frac{p-1}{2}}-1}{5}+\frac{6}{5}p\,q_p(2)\pmod{p^2}.
\]

For the analogous sums over the shorter range $0<k<p/3$ we only have the method of Section~\ref{sec:pounds} at our disposal,
which only yields congruences modulo $p$.
In particular, we can apply Equation~\eqref{eq:3kksuk-pol-short}, which follows from Theorem~\ref{thm:rkksuk_short} as explained in Remark~\ref{rem:beta}.
Thus, for $r=3$ and $x=2$ we take $\beta=(1+i)/2$, whence $c=c_1=\beta(1-\beta)=1/2$.
Because of the congruences
$\pounds_1(1/2)\equiv q_p(2)\pmod{p}$
and
$\pounds_1((1+i)/2)\equiv q_p(2)/2\pmod{p}$,
which are immediate consequences of Equation~\eqref{eq:pth_power_mod_p^2_var},
Equation~\eqref{eq:3kksuk-pol-short} yields
\[
\sum_{0<k<p/3}
\binom{3k}{k}\frac{2^k}{k}\equiv
-3q_p(2)\pmod{p}.
\]

In our next application we take $r=3$ and $x=1/8$.
Thus, we have $c_1=-1$, $c_2=-2+\sqrt{5}$, $c_3=-2-\sqrt{5}$.
Because $c_2$ and $c_3$ are harder to deal with in this case, we start with applying the results of Section~\ref{sec:pounds}.
Taking $\beta=\varphi_+$, where $\varphi_\pm=(1\pm\sqrt{5})/2$, we have $c=c_1=\beta(1-\beta)=-1$.
We have $\LL_1(-1)\equiv -2q_p(2)\pmod{p}$ and
$\LL_1(\varphi_+)\equiv(1-\varphi_+^p-\varphi_-^p)/p=-q_L\pmod{p}$.
Here the notation $q_L$ is adopted because $L_n=\varphi_+^n+\varphi_-^n$ is the $n$th Lucas number.
With these ingredients, Equation~\eqref{eq:3kksuk-pol-long} yields
\begin{align*}
\sum_{0<k<p} \binom{3k}{k} \frac{(1/8)^k}{k}
&\equiv 3(1-c)^{-2p}\bigl(\LL_1(\beta)-(1-c^p)\LL_1(c)\bigr)\pmod{p}\\
&\equiv\frac{3}{4}(\LL_1(\varphi_+)-2\LL_1(-1))
\equiv 3q_p(2)-\frac{3}{4}q_L\pmod{p},
\end{align*}
and Equation~\eqref{eq:3kksuk-pol-short} similarly yields
\[
\sum_{0<k<p/3} \binom{3k}{k} \frac{(1/8)^k}{k}
\equiv 3q_p(2)-\frac{3}{2}q_L\pmod{p}.
\]
As to the results of Section~\ref{sec:mod_p^2}, here we limit ourselves to Theorem~\ref{thm:rkkmod2}, which does not involve $\LL_2$.
We have $c_2=-\varphi_-^3$, $c_3=-\varphi_+^3$, $1-c_2=2\varphi_-^2$, $1-c_3=2\varphi_+^2$.
Also, according to~\cite[Equation~(2.9)]{Williams}, for $p\not=5$ we have
\[
L_{p-\leg{p}{5}}\equiv 2\leg{p}{5}\pmod{p^2},
\]
which implies
\[
2\varphi_{\pm}^p\equiv 1\pm\sqrt{5}\leg{p}{5}+\biggl(1\pm\frac{\sqrt{5}}{5}\leg{p}{5}\biggr)p\, q_L\pmod{p^2}.
\]
Then Theorem~\ref{thm:rkkmod2} yields
\[
\sum_{0\le k<p} \binom{3k}{k} (1/8)^k\equiv \frac{1}{4}+\frac{3}{4}\leg{p}{5}+\frac{9}{10}\leg{p}{5} p\,q_L\pmod{p^2}.
\]

Next, we set $r=3$ and $x=4/27$, whence $c_1=c_2=-2$, $c_3=1/4$.
Taking $\beta=-1$ we have $c=c_1=\beta(1-\beta)=-2$.
Since
$\pounds_1(-2)\equiv 2q_p(2)-3q_p(3)\pmod{p}$, if $p\not=3$ then
Equations~\eqref{eq:rkksuk_long} and~\eqref{eq:3kksuk-pol-short} yield
\begin{align*}
\sum_{0<k<p} \binom{3k}{k} \frac{(4/27)^k}{k}
&\equiv -\frac{8}{3}q_p(2)+3q_p(3)\pmod{p},\\
\sum_{0<k<p/3} \binom{3k}{k} \frac{(4/27)^k}{k}
&\equiv -4q_p(2)+3q_p(3)\pmod{p}.
\end{align*}
This is one case where Theorem~\ref{thm:rkkmod2_multiple} applies, with $x_0=4/27$,
and yields
\[
\sum_{0\le k<p} \binom{3k}{k}(4/27)^k
\equiv \frac{1}{9}+\frac{8}{27}p(3+q_p(2))\pmod{p^2}.
\]

As a further application of Theorem~\ref{thm:rkkmod2_multiple} we take $r=4$,
whence $x_0=27/256$.
Then we find
$c_1=c_2=-3$, $c_3=(7+4i\sqrt{2})/27$, $c_4=(7-4i\sqrt{2})/27$.
Using
$$(7\pm 4i\sqrt{2})^p\equiv 7+4i\sqrt{2}\leg{-2}{p}\pmod{p}$$
Theorem~\ref{thm:rkkmod2_multiple} yields
$$\sum_{0\le k<p} \binom{4k}{k}(27/256)^k
\equiv \frac{11}{72}+\frac{1}{288}\leg{-2}{p}\pmod{p}.$$
Here we have limited ourselves to a congruence modulo $p$
because $\LL(c_3)$ and $\LL(c_4)$ are not easy to evaluate.

Even in the well-trodden case $r=2$ our results of Section~\ref{sec:mod_p^2} yields numerical congruences that were not accessible to the methods of~\cite{MatTau:polylog}.
We list some of those, omitting the calculation details.
When $x=1/3$ we have $c_1=(-1+i\sqrt{3})/2$,
$c_2=(-1-i\sqrt{3})/2$,
and for $p>3$ we find
\begin{align*}
\sum_{0<k<p} \binom{2k}{k}\frac{(1/3)^k}{k}
&\equiv q_p(3)-\frac{1}{2}p\, q_p^2(3)\pmod{p^2},
\\
\sum_{0<k<p} \binom{2k}{k}\frac{(1/3)^k}{k^2}
&\equiv \frac{1}{9}\leg{p}{3} B_{p-2}(1/3)-\frac{1}{2}\, q_p^2(3)\pmod{p}.
\end{align*}
When $x=-2$ we have $c_1=2$, $c_2=1/2$,
and for $p>3$ we find
\begin{align*}
\sum_{0<k<p} \binom{2k}{k}\frac{(-2)^k}{k}
&\equiv -4q_p(2)+4p\, q_p^2(2)\pmod{p^2},
\\
\sum_{0<k<p} \binom{2k}{k}\frac{(-2)^k}{k^2}
&\equiv -2\, q_p^2(2)\pmod{p}.
\end{align*}
Note that the sum
$\sum_{0<k<p} \binom{2k}{k}x^k$
was already evaluated modulo $p^3$ when $x=1/3$ or $-2$ in~\cite[Section~8]{MatTau:polylog}.

\end{document}